\documentclass[11pt]{amsart}

\usepackage{amsmath}
\usepackage{amscd}
\usepackage{amssymb}
\usepackage{latexsym}
\usepackage{stmaryrd} 

\usepackage[arrow,curve,matrix,tips,2cell]{xy}
  \SelectTips{eu}{10} \UseTips
  \UseAllTwocells
\newtheorem{theorem}{Theorem}[section]
\newtheorem*{theorem*}{Theorem}
\newtheorem{lemma}[theorem]{Lemma}
\newtheorem*{lemma*}{Lemma}
\newtheorem{corollary}[theorem]{Corollary}
\newtheorem*{corollary*}{Corollary}
\newtheorem{proposition}[theorem]{Proposition}

\newtheorem{remark}[theorem]{Remark}
\newtheorem{definition}[theorem]{Definition}

\newcommand{\bgl}{\begin{equation}} 
\newcommand{\egl}{\end{equation}}
\newcommand{\bgloz}{\begin{equation*}} 
\newcommand{\egloz}{\end{equation*}}
\newcommand{\bgln}{\begin{eqnarray}} 
\newcommand{\egln}{\end{eqnarray}}
\newcommand{\bglnoz}{\begin{eqnarray*}} 
\newcommand{\eglnoz}{\end{eqnarray*}}
\newcommand{\btheo}{\begin{theorem}}
\newcommand{\etheo}{\end{theorem}}
\newcommand{\btheooz}{\begin{theorem*}}
\newcommand{\etheooz}{\end{theorem*}}
\newcommand{\blemma}{\begin{lemma}}
\newcommand{\elemma}{\end{lemma}}
\newcommand{\blemmaoz}{\begin{lemma*}}
\newcommand{\elemmaoz}{\end{lemma*}}
\newcommand{\bproof}{\begin{proof}}
\newcommand{\eproof}{\end{proof}}
\newcommand{\bbew}{\begin{beweis}}
\newcommand{\ebew}{\end{beweis}}
\newcommand{\bremark}{\begin{remark}\em}
\newcommand{\eremark}{\end{remark}}
\newcommand{\bdefin}{\begin{definition}}
\newcommand{\edefin}{\end{definition}}
\newcommand{\bprop}{\begin{proposition}}
\newcommand{\eprop}{\end{proposition}}
\newcommand{\bcor}{\begin{corollary}}
\newcommand{\ecor}{\end{corollary}}
\newcommand{\bcoroz}{\begin{corollary*}}
\newcommand{\ecoroz}{\end{corollary*}}
\newcommand{\bfa}{\begin{cases}} 
\newcommand{\efa}{\end{cases}}
\newcommand{\cC}{\mathcal C}

\newcommand{\cL}{\mathcal L}
\newcommand{\cM}{\mathcal M}

\newcommand{\cS}{\mathcal S}

\newcommand{\ba}{{\bf a}}
\newcommand{\bb}{{\bf b}}
\newcommand{\bc}{{\bf c}}
\newcommand{\bd}{{\bf d}}
\def\Az{\mathbb{A}}

\def\Nz{\mathbb{N}}

\def\Zz{\mathbb{Z}}

\def\1z{\mathbb{1}}
\newcommand{\mfa}{\mathfrak a}
\newcommand{\mfb}{\mathfrak b}
\newcommand{\mfk}{\mathfrak k}
\newcommand{\mfp}{\mathfrak p}
\newcommand{\mfq}{\mathfrak q}
\newcommand{\lori}{\longrightarrow}

\newcommand{\ma}{\mapsto} 
\newcommand\onto{\twoheadrightarrow} 
\newcommand\into{\hookrightarrow} 

\def\SEMI{\mbox{$\times\kern-2pt\vrule height5pt width.6pt \kern3pt $}}

\newcommand{\img}{{\rm Im\,}}

\newcommand{\id}{{\rm id}}

\newcommand{\coker}{{\rm coker}\,}
\newcommand{\reg}{^\times} 
\newcommand{\lspan}{{\rm span}} 
\newcommand{\clspan}{\overline{\lspan}} 
\newcommand{\Zspan}{\Zz \text{-} \lspan}
\newcommand{\abs}[1]{\lvert#1\rvert} 
\newcommand{\defeq}{\mathrel{:=}} 

\newcommand{\dop}{\text{: }} 
\newcommand{\falls}{\text{ if }} 
\newcommand{\sonst}{\text{ else}} 
\newcommand{\dotcup}{\ensuremath{\mathaccent\cdot\cup}} 

\newcommand{\lge}{\left\{} 
\newcommand{\rge}{\right\}} 
\newcommand{\lru}{\left(} 
\newcommand{\rru}{\right)} 
\newcommand{\leck}{\left[} 
\newcommand{\reck}{\right]} 
\newcommand{\lsp}{\left\langle} 
\newcommand{\rsp}{\right\rangle} 
\newcommand{\rukl}[1]{\lru #1 \rru} 
\newcommand{\eckl}[1]{\leck #1 \reck} 
\newcommand{\gekl}[1]{\lge #1 \rge} 
\newcommand{\spkl}[1]{\lsp #1 \rsp} 
\newcommand{\menge}[2]{\gekl{ #1 \dop #2 }} 
\begin{document}

\title[On K-theoretic invariants of semigroup C*-algebras, Part II]{On K-theoretic invariants of semigroup C*-algebras attached to number fields, Part II}

\author{Xin Li}

\address{Xin Li, School of Mathematical Sciences, Queen Mary University of London, Mile End Road, London E1 4NS, United Kingdom}
\email{xin.li@qmul.ac.uk}


\begin{abstract}
This paper continues the study of K-theoretic invariants for semigroup C*-algebras attached to $ax+b$-semigroups over rings of algebraic integers in number fields. We show that from the semigroup C*-algebra together with its canonical commutative subalgebra, it is possible to reconstruct the zeta function of the underlying number field as well as its ideal class group (as a group). In addition, we give an alternative interpretation of this result in terms of dynamical systems.
\end{abstract}

\thanks{Research supported by EPSRC grant EP/M009718/1.}

\maketitle


\setlength{\parindent}{0pt} \setlength{\parskip}{0.5cm}

\section{Introduction}

This paper is a sequel to \cite{Li}. We continue the study of K-theoretic invariants for semigroup C*-algebras attached to $ax+b$-semigroups over rings of algebraic integers in number fields.

The K-groups as such have been computed in \cite{C-E-L1}. The final answer involved the ideal class group of the corresponding number field. However, the class group only appeared as a set. Thus, a natural question is how to recover the group structure.

In \cite{Li}, it was shown that given the number of roots of unity, it is possible to recover the Dedekind zeta function of the underlying number field by studying the ideal structure and refined K-theoretic invariants of the corresponding semigroup C*-algebras.

In the present paper, we show that the semigroup C*-algebra together with its canonical commutative subalgebra completely determines the ideal class group. We even provide an explicit way to reconstruct the ideal class group, as a group, from the given C*-algebraic data. In addition, we can recover the zeta function without having to know the number of roots of unity.

Let us make these statements more precise and formulate our main theorem. Let $R$ be the ring of algebraic integers in a number field, and let $R \rtimes R\reg$ be the $ax+b$-semigroup over $R$. The semigroup C*-algebra $C^*(R \rtimes R\reg)$ is the concrete C*-algebra of operators on $\ell^2(R \rtimes R\reg)$ generated by (the isometries from) the left regular representation of $R \rtimes R\reg$. The algebra $\ell^{\infty}(R \rtimes R\reg)$ acts on $\ell^2(R \rtimes R\reg)$ as well, as multiplication operators, and the canonical commutative subalgebra of $C^*(R \rtimes R\reg)$ is given by $D(R \rtimes R\reg) \defeq C^*(R \rtimes R\reg) \cap \ell^{\infty}(R \rtimes R\reg)$. This intersection is taken in $\cL(\ell^2(R \rtimes R\reg))$. Now let $K$ and $L$ be two number fields with Dedekind zeta functions $\zeta_K$ and $\zeta_L$ as well as ideal class groups $Cl_K$ and $Cl_L$. Let $R$ and $S$ be the rings of algebraic integers in $K$ and $L$. We consider the semigroup C*-algebras $C^*(R \rtimes R\reg)$ and $C^*(S \rtimes S\reg)$ with their canonical commutative subalgebras $D(R \rtimes R\reg)$ and $D(S \rtimes S\reg)$, respectively. Here is our main result:
\btheo
\label{mainTHM}
If there exists an isomorphism $C^*(R \rtimes R\reg) \overset{\cong}{\lori} C^*(S \rtimes S\reg)$ which restricts to an isomorphism $D(R \rtimes R\reg) \overset{\cong}{\lori} D(S \rtimes S\reg)$, then $\zeta_K = \zeta_L$ and we have a group isomorphism $Cl_K \cong Cl_L$.
\etheo
In the proof of this result, we will explicitly reconstruct the ideal class group as a group from a C*-isomorphism preserving the commutative subalgebras.

Two number fields with the same zeta function are called arithmetically equivalent (see \cite{Per}). There are examples of arithmetically equivalent number fields which have different class numbers (see \cite{dS-P}). Therefore, the conclusion in Theorem~\ref{mainTHM} is really stronger than in \cite[Theorem~1.1]{Li}. But of course, the assumption in our present theorem is stronger as well: We ask for an isomorphism of semigroup C*-algebras which identifies the canonical commutative subalgebras. It turns out that for any ring of algebraic integers $R$, the pair $(C^*(R \rtimes R\reg), D(R \rtimes R\reg))$ is a Cartan pair, in the sense of \cite{R}. With this observation in mind, we can now give an alternative interpretation of Theorem~\ref{mainTHM} in terms of dynamical systems.

As explained in \cite{C-D-L}, given the ring of algebraic integers $R$ in a number field $K$, the semigroup C*-algebra $C^*(R \rtimes R\reg)$ is canonically isomorphic to a full corner of the crossed product $C_0(\Omega(K)) \rtimes (K \rtimes K\reg)$. Here $\Omega(K)$ is given as follows: Let $\Az_{f,K}$ be the finite adele ring over $K$, and let $\hat{R}$ be the subring of finite integral adeles. We define an equivalence relation on $\Az_{f,K} \times \Az_{f,K}$ by setting $(\bb,\ba) \sim (\bd,\bc)$ if $\ba \hat{R} = \bc \hat{R}$ and $\bb - \bd \in \ba \hat{R} = \bc \hat{R}$. $\Omega(K)$ is the quotient space $\Az_{f,K} \times \Az_{f,K} / \sim$. The action of the $ax+b$-group $K \rtimes K\reg$ on $\Az_{f,K} \times \Az_{f,K}$ given by $(b,a).(\bb,\ba) = (b + a \bb,a \ba)$ descends to an action on $\Omega(K)$. This action gives rise to the crossed product $C_0(\Omega(K)) \rtimes (K \rtimes K\reg)$ from above.

Now suppose that $G \curvearrowright X$ and $H \curvearrowright Y$ are two topological dynamical systems consisting of discrete groups $G$ and $H$ acting on locally compact Hausdorff spaces $X$ and $Y$ via homeomorphisms. We denote the actions by $G \times X \to X, \, (g,x) \ma g.x$ and $H \times Y \to Y, \, (h,y) \ma h.y$. We say that $G \curvearrowright X$ and $H \curvearrowright Y$ are continuously orbit equivalent if there exists a homeomorphism $\varphi: \: X \overset{\cong}{\lori} Y$ and continuous maps $a: \: G \times X \to H$ and $b: \: H \times Y \to G$ such that $\varphi(g.x) = a(g,x).\varphi(x)$ and $\varphi^{-1}(h.y) = \phi(h,y).\varphi^{-1}(y)$ for all $g \in G$, $h \in H$, $x \in X$ and $y \in Y$. We refer the reader to \cite{Li-pre} for more information about continuous orbit equivalence.

With these notations, here is the alternative interpretation of Theorem~\ref{mainTHM}:
\btheo
\label{THM_DS}
Let $K$ and $L$ be two number fields. If $K \rtimes K\reg \curvearrowright \Omega(K)$ and $L \rtimes L\reg \curvearrowright \Omega(L)$ are continuously orbit equivalent, then $\zeta_K = \zeta_L$ and there is a group isomorphism $Cl_K \cong Cl_L$.
\etheo

Our proofs of both theorems use C*-algebraic techniques, in particular K-theory. It would be interesting to find a number-theoretic explanation for Theorem~\ref{THM_DS}. Moreover, the following natural question is left open: Do there exist number fields $K$ and $L$ which are not isomorphic, but for which $K \rtimes K\reg \curvearrowright \Omega(K)$ and $L \rtimes L\reg \curvearrowright \Omega(L)$ are continuously orbit equivalent?

\section{K-theoretic invariants for semigroup C*-algebras and their commutative subalgebras}
\label{sec_K}

Let us briefly recall the construction of semigroup C*-algebras. Let $P$ be a left cancellative semigroup, and let $\menge{\delta_x}{x \in P}$ be the canonical orthonormal basis of $\ell^2 P$. For $p \in P$, let $V_p$ be the isometry given by $V_p: \: \ell^2 P \to \ell^2 P, \, \delta_x \ma \delta_{px}$. We define $C^*_{\lambda}(P) \defeq C^*(\menge{V_p}{p \in P}) \subseteq \cL(\ell^2 P)$, and $D_{\lambda}(P) \defeq C^*_{\lambda}(P) \cap \ell^{\infty}(P) \subseteq \cL(\ell^2 P)$.

Now let $R$ be the ring of algebraic integers in a number field $K$, and let $P = R \rtimes R\reg$ be the $ax+b$-semigroup over $R$. As a set, $R \rtimes R\reg$ is given by the direct product $R \times R\reg$, and multiplication is given by $(b,a) (d,c) = (b + ad,ac)$. We write $C^*(R \rtimes R\reg)$ for $C^*_{\lambda}(R \rtimes R\reg)$ and $D(R \rtimes R\reg)$ for $D_{\lambda}(R \rtimes R\reg)$.

Following the notation in \cite{Li}, let $E_{(r+I) \times I\reg}$, for $r \in R$ and $(0) \neq I \triangleleft R$, be the multiplication operator of the characteristic function $1_{(r+I) \times I\reg} \in \ell^{\infty}(R \rtimes R\reg)$. All these $E_{(r+I) \times I\reg}$ are elements of $D(R \rtimes R\reg)$. We abbreviate $e_{r+I} \defeq E_{(r+I) \times I\reg}$. As observed in \cite{Li}, every primitive ideal of $C^*(R \rtimes R\reg)$ is of the form
$$
  I_F \defeq \spkl{\menge{1 - \sum_{r \in R / \mfp} e_{r+\mfp}}{\mfp \in F}} \triangleleft C^*(R \rtimes R\reg)
$$
for a collection $F$ of prime ideals of $R$. Note that by convention, a prime ideal (of $R$) is always nonzero.

For the sake of brevity, set $D \defeq D(R \rtimes R\reg)$. For every subset $F$ of prime ideals of $R$, set $D_F \defeq I_F \cap D$, and let $\iota_F: \ D_F \into D$, $i: D \into C^*(R \rtimes R\reg)$ be the canonical embeddings. We denote the induced homomorphisms in $K_0$ by $(\iota_F)_*$ and $i_*$. Let $\Delta \defeq i_*(K_0(D))$, $\Delta_F \defeq i_*((\iota_F)_*(K_0(D_F)))$ and write $\pi_F$ for the canonical projection $\Delta \onto \Delta / \Delta_F$. For $F = \emptyset$, we set $D_F \defeq (0)$ and $\Delta_F = \gekl{0}$. Let $Cl$ be the class group of $K$. Given a collection $F$ of prime ideals of $R$, let $Cl_F$ be the subgroup of $Cl$ given by $\spkl{\menge{[\mfp]}{\mfp \in F}} \subseteq Cl$, where $Cl_{\emptyset}$ is the trivial subgroup. The following is the crucial technical result of this paper:
\bprop
\label{D/D?}
Let $F$ be a finite collection of prime ideals $\mfp$ of $R$ with $\mfp \nmid 2$.

There exists $d_F \in \Nz_0$ with $\Delta / \Delta_F \cong \bigoplus_{Cl / Cl_F} \Zz / d_F \Zz$. For $F = \emptyset$, we have $d_{\emptyset} = 0$, and for $F \neq \emptyset$, $d_F$ is positive and even.

Moreover, there are ideals $\mfa_{F,i}$ of $R$ with $Cl = \dotcup_i \, Cl_F [\mfa_{F,i}]$ and such that $\gekl{\pi_F [e_{\mfa_{F,i}}]}_i$ forms a $(\Zz / d_F \Zz)$-basis for $\Delta / \Delta_F$. Given an ideal $\mfa$ of $R$ with $[\mfa] \in Cl_F [\mfa_{F,i}]$, there exists an odd number $l_F(a) \in \Nz$ with $\pi_F [e_{\mfa}] = l_F(\mfa) \pi_F [e_{\mfa_{F,i}}]$.
\eprop
We will prove this proposition by induction on $\abs{F}$. For the induction step, we need a bit of preparation. From \cite{C-E-L1} (see also \cite{Li}), we obtain a canonical identification $\Delta \cong \Zz[Cl]$ as abelian groups. Every prime ideal $\mfp$ induces a homomorphism $\mfp_*: \ \Delta \cong \Zz[Cl] \to \Zz[Cl] \cong \Delta$, $[e_{\mfa}] \ma [e_{\mfp \mfa}]$. Let $M_{\mfp}$ be the homomorphism $\Delta \overset{\mfp_*}{\lori} \Delta \overset{N(\mfp) \, \id}{\lori} \Delta$, $[e_{\mfa}] \ma N(\mfp) [e_{\mfp \mfa}]$.
\blemma
\label{D_F=SUM}
$\Delta_F = \sum_{\mfp \in F} (\id - M_{\mfp})(\Delta)$.
\elemma
\bproof
By \cite{E-L} and the construction of $D_F$, we know that
$$D_F = C^*(\menge{e_{x+\mfa} - \sum_{r \in \mfa / \mfp \mfa} e_{x + r + \mfp \mfa}}{x \in R, \, (0) \neq \mfa \triangleleft R, \, \mfp \in F}) \subseteq D.$$
Since $\lspan(\menge{e_{x+\mfa} - \sum_{r \in \mfa / \mfp \mfa} e_{x + r + \mfp \mfa}}{x \in R, \, (0) \neq \mfa \triangleleft R, \, \mfp \in F})$ is multiplicatively closed, we conclude that
$$D_F = \clspan(\menge{e_{x+\mfa} - \sum_{r \in \mfa / \mfp \mfa} e_{x + r + \mfp \mfa}}{x \in R, \, (0) \neq \mfa \triangleleft R, \, \mfp \in F}).$$
Let $\cS_F \defeq \menge{e_{x+\mfa} - \sum_{r \in \mfa / \mfp \mfa} e_{x + r + \mfp \mfa}}{x \in R, \, (0) \neq \mfa \triangleleft R, \, \mfp \in F}$ and set $[\cS_F] \defeq \menge{[a]}{a \in \cS_F} \subseteq K_0(D)$. Moreover, let $\cM_F$ be the multiplicatively closed set of projections in $D_F$ generated by $\cS_F$, and set $[\cM_F] \defeq \menge{[a]}{a \in \cM_F} \subseteq K_0(D)$. By \cite{C-E-L1}, we have $(\iota_F)_*(K_0(D_F)) = \Zspan([\cM_F])$. To show that $\Zspan([\cM_F]) = \Zspan([\cS_F])$, we show that for all $a$ and $b$ in $\cS_F$, we have $[ab] \in \Zspan([\cS_F])$. Let $a = e_{x+\mfa} - \sum_r e_{x + r + \mfp \mfa}$, $b = e_{y+\mfb} - \sum_s e_{y + s + \mfq \mfb}$. Without loss of generality, we may assume $e_{x+\mfa} e_{y+\mfb} \neq 0$, i.e., $(x+\mfa) \cap (y+\mfb) \neq \emptyset$. Replacing $x$ and $y$ by $z \in (x+\mfa) \cap (y+\mfb)$, we may assume without loss of generality that $x=y$. And since $\cS_F$ is invariant under the additive action of $R$, we may further assume $x=y=0$. So $a = e_{\mfa} - \sum_r e_{r + \mfp \mfa}$ and $b = e_{\mfb} - \sum_s e_{s + \mfq \mfb}$. As $ab = 0$ if $v_{\mfp}(\mfb) > v_{\mfp}(\mfa)$ or $v_{\mfq}(\mfa) > v_{\mfq}(\mfb)$, we may assume $v_{\mfp}(\mfb) \leq v_{\mfp}(\mfa)$ and $v_{\mfq}(\mfa) \leq v_{\mfq}(\mfb)$. Then $(\mfp \mfa) \cap \mfb = \mfp (\mfa \cap \mfb)$, $\mfa \cap (\mfq \mfb) = \mfq (\mfa \cap \mfb)$ and 
$(\mfp \mfa) \cap (\mfq \mfb) =
  \bfa
    \mfp \mfq (\mfa \cap \mfb) \falls \mfp \neq \mfq \\
    \mfp (\mfa \cap \mfb) \falls \mfp = \mfq
  \efa.
$
Therefore, we have in case $\mfp \neq \mfq$:
\bgloz
  ab = \rukl{e_{\mfa \cap \mfb} - \sum_{\rho} e_{\rho + \mfp (\mfa \cap \mfb)}}
  - \sum_{\sigma} \rukl{e_{\sigma + \mfq (\mfa \cap \mfb)} - \sum_{\tau_{\sigma}} e_{\tau_{\sigma} + \mfp \mfq (\mfa \cap \mfb)}}
\egloz
and for all $\sigma$,
\bgloz
  e_{\sigma + \mfq (\mfa \cap \mfb)} - \sum_{\tau_{\sigma}} e_{\tau_{\sigma} + \mfp \mfq (\mfa \cap \mfb)} 
  = \rukl{e_{\mfa \cap \mfb} - \sum_{\rho} e_{\rho + \mfp (\mfa \cap \mfb)}} \cdot e_{\sigma + \mfq (\mfa \cap \mfb)}
  \leq e_{\mfa \cap \mfb} - \sum_{\rho} e_{\rho + \mfp (\mfa \cap \mfb)}.
\egloz
Thus
\bgloz
  [ab] = \eckl{e_{\mfa \cap \mfb} - \sum_{\rho} e_{\rho + \mfp (\mfa \cap \mfb)}}
  - \sum_{\sigma} \eckl{e_{\sigma + \mfq (\mfa \cap \mfb)} - \sum_{\tau_{\sigma}} e_{\tau_{\sigma} + \mfp \mfq (\mfa \cap \mfb)}}
\egloz
lies in $\Zspan([\cS_F])$.

If $\mfp = \mfq$, then $ab = e_{\mfa \cap \mfb} - \sum_{\rho} e_{\rho + \mfp (\mfa \cap \mfb)}$. Thus $[ab] = \eckl{e_{\mfa \cap \mfb} - \sum_{\rho} e_{\rho + \mfp (\mfa \cap \mfb)}}$ lies in $\Zspan([\cS_F])$.

Since $\eckl{e_{x+\mfa} - \sum_{r \in \mfa / \mfp \mfa} e_{x + r + \mfp \mfa}} = [e_{\mfa}] - N(\mfp) [e_{\mfp \mfa}]$ holds in $K_0(D)$, our claim follows.
\eproof

\blemma
\label{coker}
Let $\alpha$ be an endomorphism of the finitely generated abelian group $\bigoplus_{i=1}^n \Zz / d \Zz$ ($d \in \Nz_0$). Assume that with respect to the standard $\Zz / d \Zz$-basis $e_i$, $\alpha$ is given by the matrix
$
  \rukl{
  \begin{smallmatrix}
  1 & 0 & \dotso & 0 & -N_n \\
  -N_1 & 1 & & 0 & 0 \\
  0 & -N_2 & & & \vdots \\
  \vdots & & \ddots & & \vdots \\
  \vdots & & & 1 & 0 \\
  0 & \dotso & \dotso & -N_{n-1} & 1
  \end{smallmatrix}
  }
$.
Then with $d_{\alpha} \defeq \gcd(N_1 \dotsm N_n - 1, d)$, we have $\coker(\alpha) \cong \Zz / d_{\alpha} \Zz$. Moreover, if $\pi: \ \bigoplus_{i=1}^n \Zz / d \Zz \onto \coker(\alpha)$ is the canonical projection, then $\pi(e_n)$ is a generator of $\coker(\alpha)$. For every $1 \leq i \leq n$, we have $\pi(e_i) = \rukl{\prod_{h=1}^i N_h} \pi(e_n)$.
\elemma
\bproof
Let $f_i \defeq e_i - N_i e_{i+1}$ for $1 \leq i \leq n-1$ and $f_n \defeq e_n$. Obviously, $\gekl{f_i}_{1 \leq i \leq n}$ forms a $\Zz / d \Zz$-basis of $\bigoplus_{i=1}^n \Zz / d \Zz$. We have $\img(\alpha) = \sum_{i=1}^{n-1} \Zz f_i + (N_1 \dotsm N_n - 1) \Zz f_n$ since we obviously have $f_i \in \img(\alpha)$, and if $\alpha(\sum_i m_i e_i) = m e_n$, then
\bgloz
  \sum_{i=1}^{n-1} (m_i e_i - N_i m_i e_{i+1}) + m_n e_n - N_n m_n e_1 = m e_n,
\egloz
thus
\bgloz
  m_i - N_{i-1} m_{i-1} \equiv 0 \mod d \Zz \ {\rm for} \ {\rm all} \ 2 \leq i \leq n-1, \ {\rm and} \ m_1 - N_n m_n \equiv 0 \mod d \Zz.
\egloz
Therefore,
\bglnoz
  && m_{n-1} \equiv N_{n-2} m_{n-2} \equiv \dotso \equiv (\prod_{i \neq n-1} N_i) \cdot m_n \mod d \Zz \\
{\rm and} && m \equiv m_n - N_{n-1} m_{n-1} \equiv (1 - \prod_{i=1}^n N_i) \cdot m_n \mod d \Zz.
\eglnoz
\eproof
After these preparations, we are ready for the proof of our proposition. 
\bproof[Proof of Proposition~\ref{D/D?}]
As we mentioned obove, we proceed inductively on $\abs{F}$. The case $F = \emptyset$ is just the canonical identification $\Delta \cong \Zz[Cl]$ from \cite{C-E-L1} (see also \cite{Li}). Now let us assume that we know our assumption for $F_m = \gekl{\mfp_1, \dotsc, \mfp_m}$. Set $\Delta_m \defeq \Delta_{F_m}$. For another prime ideal $\mfp_{m+1}$ (with $\mfp_{m+1} \nmid 2$), set $F_{m+1} = \gekl{\mfp_1, \dotsc, \mfp_{m+1}}$ and $\Delta_{m+1} \defeq \Delta_{F_{m+1}}$. Since $M_{m+1} \defeq M_{\mfp_{m+1}}$ commutes with $M_{\mfp}$ for every $\mfp \in \gekl{\mfp_1, \dotsc, \mfp_m}$, $M_{m+1}$ induces a homomorphism $M_{m+1}^{\bullet}: \ \Delta / \Delta_m \to \Delta / \Delta_m$. It is then obvious from Lemma~\ref{D_F=SUM} that there is a canonical isomorphism $\Delta / \Delta_{m+1} \cong \coker(\id_{\Delta / \Delta_m} - M_{m+1}^{\bullet})$. Let $Cl_m \defeq Cl_{F_m}$, $Cl_{m+1} \defeq Cl_{F_{m+1}}$, $d_m \defeq d_{F_m}$, $d_{m+1} \defeq d_{F_{m+1}}$ and $\pi_m \defeq \pi_{F_m}$. Moreover, we choose ideals $\mfa_{F_m,i}$ as in our proposition for $F_m$ and set $\mfa_{m,i} \defeq \mfa_{F_m,i}$. Let $\gekl{\mfa_{m,i}} = \gekl{\mfa_{m,j,k}}$ ($j=0,1,\dotsc$ and $k=0,1,\dotsc$) such that $Cl = \dotcup_j Cl_{m+1} [\mfa_{m,j,0}]$ and $Cl_{m+1} [\mfa_{m,j,0}] = \dotcup_k Cl_m [\mfa_{m,j,k}]$ with $Cl_m [\mfa_{m,j,k}] = Cl_m [\mfp_{m+1}]^k [\mfa_{m,j,0}]$. Then $\id_{\Delta / \Delta_m} - M_{m+1}^{\bullet}$ respects the direct sum decomposition
$$
  \Delta / \Delta_m = \bigoplus_i (\Zz / d_m \Zz) \pi_m [e_{\mfa_{m,i}}] = \bigoplus_j (\bigoplus_k (\Zz / d_m \Zz) \pi_m [e_{\mfa_{m,j,k}}]).
$$
So
\bgl
\label{coker=oplus}
  \coker(\id_{\Delta / \Delta_m} - M_{m+1}^{\bullet})
  \cong \bigoplus_j \coker \rukl{(\id_{\Delta / \Delta_m} - M_{m+1}^{\bullet}) \vert_{\bigoplus_k (\Zz / d_m \Zz) \pi_m [e_{\mfa_{m,j,k}}]}}
\egl
Moreover, for arbitrary $j \neq j'$, the automorphism $Cl \to Cl$, $\mfk \ma [\mfa_{m,j',0}][\mfa_{m,j,0}]^{-1} \mfk$ induces an automorphism of $\Delta$ sending $\bigoplus_k \Zz [e_{\mfa_{m,j,k}}]$ to $\bigoplus_k \Zz [e_{\mfa_{m,j',k}}]$. This automorphism commutes with $\id - M_{\mfp}$ for every $\mfp \in \gekl{\mfp_1, \dotsc, \mfp_m}$, so it induces an automorphism of $\Delta / \Delta_m$ sending $\bigoplus_k (\Zz / d_m \Zz) [e_{\mfa_{m,j,k}}]$ to $\bigoplus_k (\Zz / d_m \Zz) [e_{\mfa_{m,j',k}}]$. As this induced automorphism also commutes with $\id_{\Delta / \Delta_m} - M_{m+1}^{\bullet}$, it in turn induces an automorphism of $\Delta / \Delta_{m+1} \cong \coker(\id_{\Delta / \Delta_m} - M_{m+1}^{\bullet})$ sending
$$\coker \rukl{(\id_{\Delta / \Delta_m} - M_{m+1}^{\bullet}) \vert_{\bigoplus_k (\Zz / d_m \Zz) \pi_m [e_{\mfa_{m,j,k}}]}}$$
to
$$\coker \rukl{(\id_{\Delta / \Delta_m} - M_{m+1}^{\bullet}) \vert_{\bigoplus_k (\Zz / d_m \Zz) \pi_m [e_{\mfa_{m,j',k}}]}}.$$
Thus all the individual summands in \eqref{coker=oplus} are isomorphic. Therefore, it suffices to determine one of them.

Let us now fix $j$. For an ideal $\mfa$ with $[\mfa] \in Cl_m [\mfa_{m,i}]$, we set $l_m(\mfa) \defeq l_{F_m}(\mfa)$, so that $\pi_m [e_{\mfa}] = l_m(\mfa) \pi_m [e_{\mfa_{m,i}}]$. Then
\bglnoz
  && M_{m+1}^{\bullet} (\pi [e_{\mfa_{m,j,k}}]) \\
  &=&
  \bfa
    N(\mfp_{m+1}) l_m(\mfp_{m+1} \mfa_{m,j,k}) \pi [e_{\mfa_{m,j,k+1}}] & \falls Cl_m [\mfp_{m+1}] [\mfa_{m,j,k}] \neq Cl_m [\mfa_{m,j,0}]; \\
    N(\mfp_{m+1}) l_m(\mfp_{m+1} \mfa_{m,j,k}) \pi [e_{\mfa_{m,j,0}}] & \sonst.
  \efa
\eglnoz
Therefore, we can apply Lemma~\ref{coker} to $\alpha = (\id_{\Delta / \Delta_m} - M_{m+1}^{\bullet}) \vert_{\bigoplus_k (\Zz / d_m \Zz) \pi_m [e_{\mfa_{m,j,k}}]}$, $d = d_m$, the $\Zz / d \Zz$-basis $e_k = \pi_m [e_{\mfa_{m,j,k}}]$ of $\bigoplus_k (\Zz / d_m \Zz) \pi_m [e_{\mfa_{m,j,k}}]$, and $N_k = N(\mfp_{m+1}) l_m(\mfp_{m+1} \mfa_{m,j,k})$. This concludes the proof of our proposition.
\eproof
A completely analogous argument yields
\bcor
\label{F=p}
$\Delta / \Delta_{\gekl{\mfp}} \cong \bigoplus_{Cl / \spkl{[\mfp]}} \Zz / (N(\mfp)^{\# \spkl{[\mfp]}} - 1) \Zz$ for all prime ideals $\mfp$.
\ecor
Note that here, $\mfp$ is an arbitrary prime ideal (not necessarily with the property $\mfp \nmid 2$).

\section{Reconstruction of norms and class groups}

We keep the notation from the previous section. The following lemma explains why it is enough to consider the prime ideals $\mfp$ of $R$ with $\mfp \nmid 2$, as we did in Proposition~\ref{D/D?}.
\blemma
\label{pnmid2}
$\menge{[\mfp]}{(0) \neq \mfp \triangleleft R \ prime \ ideal, \ \mfp \nmid 2}$ generates $Cl$.
\elemma
\bproof
Let $\gekl{\mfq_i} = \menge{\mfp}{(0) \neq \mfp \triangleleft R \ prime \ ideal, \ \mfp \mid 2}$, and choose $\mfq \in \gekl{\mfq_i}$. By Strong Approximation, there exists $x$ in $K\reg$ such that $v_{\mfq}(x) = 1$ and $v_{\mfq_i}(x) = 0$ for all $\mfq_i \neq \mfq$. Then $xR = \mfq \cdot \prod_{\mfp \nmid 2} \mfp^{v_{\mfp}(x)}$. Therefore, $[\mfq] \cdot \prod_{\mfp \nmid 2} [\mfp]^{v_{\mfp}(x)} = 1$ in $Cl$. Hence $[\mfq] \in \spkl{\menge{[\mfp]}{(0) \neq \mfp \triangleleft R \ prime \ ideal, \ \mfp \nmid 2}}$.
\eproof

Given a natural number $n \geq 1$ with prime factorization $n = \prod_p p^{v_p}$, we write $(n)_p = p^{v_p}$, or simply $n_p = p^{v_p}$. Moreover, given a subgroup $H$ of a group $G$, we write $[G:H]$ for the cardinality of the set of left cosets $G/H$. The following purely group-theoretic fact will allow us to reconstruct the class group.
\blemma
\label{C_p}
Let $C$ be a finite abelian group. Then $C$ admits a decomposition $C = \bigoplus_p C_p$ where the $C_p$ are $p$-groups, for $p$ prime. Let $\cC$ be a (possibly infinite) family of generators for $C$. For every prime $p$, we define numbers $d_p^1, d_p^2, \dotsc$ and elements $c_1, c_2, \dotsc$ of $\cC$ recursively as follows:
$$
  d_p^1 \defeq \max \rukl{\menge{(\# \spkl{c})_p}{c \in \cC}},
$$
and we choose $c_1 \in \cC$ with $(\# \spkl{c})_p = d_p^1$. For $i > 1$, we set
$$
  d_p^i \defeq \max \rukl{\menge{[\spkl{c_1, \dotsc, c_{i-1}, c} : \spkl{c_1, \dotsc, c_{i-1}}]_p}{c \in \cC \setminus \gekl{c_1, \dotsc, c_{i-1}}}},
$$
and we choose $c_i \in \cC \setminus \gekl{c_1, \dotsc, c_{i-1}}$ with $[\spkl{c_1, \dotsc, c_{i-1}, c} : \spkl{c_1, \dotsc, c_{i-1}}]_p = d_p^i$.

Then $C_p \cong \bigoplus_i \Zz / d_p^i \Zz$.
\elemma
\bproof
There exists a decreasing sequence of natural numbers $n_1 \geq n_2 \geq \dotso \geq n_j$ such that $C_p \cong \bigoplus_{i=1}^j \Zz / p^{n_i} \Zz$. We proceed inductively on $j$. The induction starts with $j=1$. In this case, our claim is obviously true. To go from $j$ to $j+1$, choose $c_1 \in \cC$ such that $(\# \spkl{c_1})_p = p^{n_1}$. Such an element always exists because $\cC$ generates $C$. Now set $C' \defeq C / \spkl{c_1}$, $\cC' \defeq \menge{c^{\cdot} \in C'}{c \in \cC \setminus \gekl{c_1}}$. Obviously, $C'_p \cong \bigoplus_{i=2}^j \Zz / p^{n_i} \Zz$. Now we can just apply the induction hypothesis to $C'$, using that $\spkl{c_1, \dotsc, c_{i-1}, c} / \spkl{c_1, \dotsc, c_{i-1}} \cong \spkl{c_2^{\cdot}, \dotsc, c_{i-1}^{\cdot}, c^{\cdot}}_{C'} / \spkl{c_2^{\cdot}, \dotsc, c_{i-1}^{\cdot}}_{C'}$ holds for all $c_2, \dotsc, c_{i-1} \in \cC$ and $c \in \cC \setminus \gekl{c_1, \dotsc, c_{i-1}}$.
\eproof

We are now ready for the proof of our main result.
\btheooz[Theorem~\ref{mainTHM}]
Let $K$ and $L$ be two number fields with rings of algebraic integers $R$ and $S$. If there exists an isomorphism $C^*(R \rtimes R\reg) \cong C^*(S \rtimes S\reg)$ which identifies $D(R \rtimes R\reg)$ with $D(S \rtimes S\reg)$, then $K$ and $L$ are arithmetically equivalent, and their class groups are isomorphic (as abelian groups), i.e., $Cl_K \cong Cl_L$.
\etheooz
\bproof
We just have to show how to read off $\menge{N(\mfp)}{(0) \neq \mfp \triangleleft R \ prime \ ideal}$ and $Cl_K$ from K-theoretic data attached to the pair $D(R \rtimes R\reg) \subseteq C^*(R \rtimes R\reg)$.

For the norms, consider the minimal nonzero primitive ideals of $C^*(R \rtimes R\reg)$. By \cite[Theorem~4.1]{Li}, these are of the form $I_{\gekl{\mfp}}$, $(0) \neq \mfp \triangleleft R$ prime ideal. Consider the canonical inclusions
$$
  I_{\gekl{\mfp}} \cap D(R \rtimes R\reg) \overset{\iota_{\gekl{\mfp}}}{\into} D(R \rtimes R\reg) \overset{i}{\into} C^*(R \rtimes R\reg)
$$
and the induced homomorphisms $(\iota_{\gekl{\mfp}})_*$ and $i_*$ in $K_0$. The quotient $\img(i_*) / \img(i_* \circ (\iota_{\gekl{\mfp}})_*)$ is isomorphic to $\bigoplus_{Cl_K / \spkl{[\mfp]}} \Zz / (N(\mfp)^{\# \spkl{[\mfp]}} - 1) \Zz$ by Corollary~\ref{F=p}. So we can extract the invariants $[Cl_K : \spkl{[\mfp]}]$ and $N(\mfp)^{\# \spkl{[\mfp]}}$. But we can also extract $\# Cl_K$ since $\img(i_*)$ has rank $\# Cl_K$ by Proposition~\ref{D/D?}. Thus we can extract $N(\mfp)$. This settles the first claim about arithmetic equivalence in our theorem (compare \cite{P-S}).

At the same time, we see that by looking at K-theory, we can distinguish between $I_{\gekl{\mfp}}$ for $\mfp \mid 2$ and $I_{\gekl{\mfp}}$ for $\mfp \nmid 2$ by checking whether $2$ divides $N(\mfp)$ or not. Using this observation, we can now read off the class group $Cl_K$ from K-theoretic invariants attached to $D(R \rtimes R\reg) \subseteq C^*(R \rtimes R\reg)$. Consider all the nonzero primitive ideals of $C^*(R \rtimes R\reg)$ which only contain finitely many nonzero primitive ideals but do not contain an ideal of the form $I_{\gekl{\mfp}}$ for $\mfp \mid 2$. These ideals are of the form $I_F$, where $F$ is a non-empty, finite subset of $\menge{\mfp}{(0) \neq \mfp \triangleleft R \ prime \ ideal, \ \mfp \nmid 2}$. This follows from \cite{E-L} using analogous arguments as in \cite[\S~3]{Li}. Again, consider the canonical inclusions
$$
  I_F \cap D(R \rtimes R\reg) \overset{\iota_F}{\into} D(R \rtimes R\reg) \overset{i}{\into} C^*(R \rtimes R\reg)
$$
and the induced homomorphisms $(\iota_F)_*$ and $i_*$ in $K_0$. By Proposition~\ref{D/D?}, the quotient $\img(i_*) / \img(i_* \circ (\iota_F)_*)$ is isomorphic to $\bigoplus_{Cl_K / \spkl{\menge{[\mfp]}{\mfp \in F}}} \Zz / d_F \Zz$. Thus from this quotient, we can extract $[Cl_K : \spkl{\menge{[\mfp]}{\mfp \in F}}]$, hence also $\# \spkl{\menge{[\mfp]}{\mfp \in F}}$. Now we can apply Lemma~\ref{C_p} to
$$
  C = Cl_K \ \text{and} \ \cC = \menge{[\mfp]}{(0) \neq \mfp \triangleleft R \ prime \ ideal, \ \mfp \nmid 2}
$$
($\cC$ generates $Cl_K$ by Lemma~\ref{pnmid2}). As a conclusion, we can extract (the isomorphism type of) $Cl_K$ from K-theoretic invariants attached to $D(R \rtimes R\reg) \subseteq C^*(R \rtimes R\reg)$, as claimed.
\eproof

\bremark
Compared with the proof of the main result in \cite{Li}, the proof of our main theorem is much cleaner: We do not need to know the number of roots of unity in advance, and we are able to extract all the norms, not only those up to a finite set of exceptions as in \cite{Li}.
\eremark

\section{Alternative interpretation in terms of dynamical systems}

In the following, we give an alternative interpretation of Theorem~\ref{mainTHM} in terms of dynamical systems which are constructed as follows: Given a number field $K$ with ring of algebraic integers $R$, let $\Az_{f,K}$ be the finite adele ring of $K$, and let $\hat{R}$ be the subring of finite integral adeles in $\Az_{f,K}$. We define an equivalence relation on $\Az_{f,K} \times \Az_{f,K}$ by setting $(\bb,\ba) \sim (\bd,\bc)$ if $\ba \hat{R} = \bc \hat{R}$ and $\bb - \bd \in \ba \hat{R} = \bc \hat{R}$. $\Omega(K)$ is the quotient space $\Az_{f,K} \times \Az_{f,K} / \sim$. We denote the class of $(\bb,\ba) \in \Az_{f,K} \times \Az_{f,K}$ in $\Az_{f,K} \times \Az_{f,K} / \sim$ by $[\bb,\ba]$. The $ax+b$-group $K \rtimes K\reg$ over $K$ acts on $\Omega(K) = \Az_{f,K} \times \Az_{f,K} / \sim$ via $(b,a).[\bb,\ba] = [b + a \bb,a \ba]$. Here we use addition and multiplication in $\Az_{f,K}$ and view $K$ and $\hat{R}$ as subrings of $\Az_{f,K}$ in the canonical way. This dynamical system $K \rtimes K\reg \curvearrowright \Omega(K)$ gives rise to the crossed product $C_0(\Omega(K)) \rtimes (K \rtimes K\reg)$.

Let $\Omega(R)$ be the set $\menge{[\bb,\ba] \in \Omega(K)}{\bb,\ba \in \hat{R}}$. It is easy to see that $\Omega(R)$ is a clopen subspace of $\Omega(K)$. As observed in \cite[\S~5]{C-D-L}, we have a canonical isomorphism $C^*(R \rtimes R\reg) \cong 1_{\Omega(R)} \rukl{C_0(\Omega(K)) \rtimes (K \rtimes K\reg)} 1_{\Omega(R)}$ sending $V_{(b,a)}$ to $1_{\Omega(R)} U_{(b,a)} 1_{\Omega(R)}$. Here we denote by the $U$ the canonical unitaries in the multiplier algebra of the crossed product implementing the (dual) action $K \rtimes K\reg \curvearrowright C_0(\Omega(K))$.

Moreover, recall the notion of continuous orbit equivalence from the introduction: Let $G \curvearrowright X$ and $H \curvearrowright Y$ be two topological dynamical systems consisting of discrete groups $G$ and $H$ acting on locally compact Hausdorff spaces $X$ and $Y$ via homeomorphisms. We denote the actions by $G \times X \to X, \, (g,x) \ma g.x$ and $H \times Y \to Y, \, (h,y) \ma h.y$. Our systems $G \curvearrowright X$ and $H \curvearrowright Y$ are called continuously orbit equivalent if there exists a homeomorphism $\varphi: \: X \overset{\cong}{\lori} Y$ and continuous maps $a: \: G \times X \to H$, $b: \: H \times Y \to G$ such that $\varphi(g.x) = a(g,x).\varphi(x)$ and $\varphi^{-1}(h.y) = b(h,y).\varphi^{-1}(y)$ for all $g \in G$, $h \in H$, $x \in X$ and $y \in Y$.

Our goal now is to prove
\btheooz[Theorem~\ref{THM_DS}]
Let $K$ and $L$ be two number fields. If $K \rtimes K\reg \curvearrowright \Omega(K)$ and $L \rtimes L\reg \curvearrowright \Omega(L)$ are continuously orbit equivalent, then $\zeta_K = \zeta_L$ and there is a group isomorphism $Cl_K \cong Cl_L$.
\etheooz
\bproof
The actions $K \rtimes K\reg \curvearrowright \Omega(K)$ and $L \rtimes L\reg \curvearrowright \Omega(L)$ are topologically free by \cite{E-L}. Therefore, \cite[Theorem~1.2]{Li-pre} tells us that $K \rtimes K\reg \curvearrowright \Omega(K)$ and $L \rtimes L\reg \curvearrowright \Omega(L)$ are continuously orbit equivalent if and only if there exists an isomorphism $C_0(\Omega(K)) \rtimes (K \rtimes K\reg) \overset{\cong}{\lori} C_0(\Omega(L)) \rtimes (L \rtimes L\reg)$ which restricts to an isomorphism $C_0(\Omega(K)) \cong C_0(\Omega(L))$.

Now, \cite{E-L} tells us that every nonzero primitive ideal of $C_0(\Omega(K)) \rtimes (K \rtimes K\reg)$ is of the form
$$
  \bar{I}_F \defeq \spkl{\menge{1_{\Omega(R)} - \sum_{r \in R / \mfp} 1_{r+\pi_{\mfp} \Omega(R)}}{\mfp \in F}} \triangleleft C_0(\Omega(K)) \rtimes (K \rtimes K\reg),
$$
where $F$ is a collection of prime ideals of $R$ and $\pi_{\mfp}$ is a finite adele with $v_{\mfp}(\pi_{\mfp}) = 1$ and $v_{\mfq}(\pi_{\mfp}) = 0$ for all $\mfq \neq \mfp$ ($v_{\mfq}$ is the valuation corresponding to $\mfq$). The connection between $\bar{I}_F$ and the ideals $I_F$ introduced in \S~\ref{sec_K} is given by the observation that $I_F \cong 1_{\Omega(R)} \bar{I}_F 1_{\Omega(R)}$ under the canonical isomorphism $C^*(R \rtimes R\reg) \cong 1_{\Omega(R)} \rukl{C_0(\Omega(K)) \rtimes (K \rtimes K\reg)} 1_{\Omega(R)}$.

This canonical isomorphism gives rise to a canonical embedding $j: \: C^*(R \rtimes R\reg) \into C_0(\Omega(K)) \rtimes (K \rtimes K\reg)$. Let $j_F: \: I_F \to \bar{I}_F$ be the restriction of $j$ to $I_F$. Moreover, we introduce the following notation, as before: Set $\bar{D} \defeq C_0(\Omega(K))$. For every subset $F$ of prime ideals of $R$, set $\bar{D}_F \defeq \bar{I}_F \cap \bar{D}$, and let $\bar{\iota}_F: \ \bar{D}_F \into \bar{D}$, $\bar{i}: \bar{D} \into C_0(\Omega(K)) \rtimes (K \rtimes K\reg)$ be the canonical embeddings. We denote the induced homomorphisms in $K_0$ by $(\bar{\iota}_F)_*$ and $\bar{i}_*$. Let $\bar{\Delta} \defeq \bar{i}_*(K_0(\bar{D}))$, $\bar{\Delta}_F \defeq \bar{i}_*((\bar{\iota}_F)_*(K_0(\bar{D}_F)))$. For $F = \emptyset$, we set $\bar{D}_F \defeq (0)$ and $\bar{\Delta}_F = \gekl{0}$.

Obviously, we have $j \circ i = \bar{i} \circ (j \vert_D)$. As $\img(\bar{i}_*) = \img(\bar{i}_* \circ (j \vert_D)_*)$ follows from \cite{C-E-L1}, we infer that $j_*(\Delta) = \bar{\Delta}$. Moreover, we also have $j \circ i \circ \iota_F = \bar{i} \circ \bar{\iota}_F \circ j_F$. It follows from \cite{C-E-L1} and an analogous argument as in the proof of Lemma~\ref{D_F=SUM} that $\img(\bar{i}_* \circ (\bar{\iota}_F)_* \circ (j_F)_*) = \img(\bar{i}_* \circ (\bar{\iota}_F)_*)$. Therefore, $j_*(\Delta_F) = \bar{\Delta}_F$. Since $j$ induces an isomorphism in $K_0$, being the embedding of a full corner, we conclude that $j_*$ gives rise to an isomorphism $\Delta / \Delta_F \cong \bar{\Delta} / \bar{\Delta}_F$, for every collection $F$ of prime ideals of $R$. 

With the help of this observation, we can now extract the zeta function and the class group of our number field $K$ from K-theoretic invariants of the pair $(C_0(\Omega(K)) \rtimes (K \rtimes K\reg), C_0(\Omega(K)))$ in exactly the same way as in the case of the pair $(C^*(R \rtimes R\reg), D(R \rtimes R\reg))$ in the proof of Theorem~\ref{mainTHM}.
\eproof

\end{document}